# Formulas for Hitting Times and Cover Times for Random Walks on Groups

Christopher Zhang

Advisor: Ramon van Handel

May 9, 2017

**Abstract**

Using the results of Ding, Lee, Peres [3], we develop formulas to compute the hitting times and cover times for random walks on groups. We developed an explicit formula for hitting times in terms of the irreducible representations of the group. We also have a way of computing cover times in terms of these hitting times. This computation is based on a quantity we indentified, which we call the volume growth function. And we believe that it is the right object to study in order to understand the cover time.

## 1 Introduction

The hitting time of a Markov chain is the time for a Markov chain to reach a certain state. These values have long been studied such as in problems like gambler's ruin. More recently, there has been interest in the cover time i.e. the time it takes for a finite Markov chain to hit all of its states. Computing cover times has applications in Markov chain algorithms such as page rank. Here we are interested in the asymptotics of the cover times of a family Markov chains, where the size of the state space goes to infinity.

Matthews [6] shows a relationship between cover times and a natural object of study, the maximum hitting time. For a Markov chain with state space $\Omega$, the Matthews upper and lower bounds differ by a factor of $\log |\Omega|$. Both the upper and lower Matthews bounds can be achieved for example with the cycle and the hypercube (see examples 4.14 and 4.16).

The paper of Ding, Lee, and Peres [3] was a big breakthrough for the theory showing the cover times for random walks on graphs are asymptotic to the expected maximum of an associated Gaussian process. Their results can be sharper than Matthews bounds as they encapsulate the information of all of the hitting times instead of only the maximum hitting time. However in the general case, computing the maximum of this Gaussian process is difficult to compute. For example, it can be computed with Talagrand's theory of majorizing measures [4].

We consider a case in which there are many symmetries, a random walk on a group. This allows us to do explicit computations of hitting times and cover times. We derive an explicit formula for hitting times for random walks on group in terms of its irreducible representations. The strength of this method is that it is exact and general. The formulas hold for any random walk on a group while conventional hitting time calculations only hold for random walks where the transition probabilities are very carefully chosen. However, the drawback is that the formulas are often difficult to calculate or bound by hand, but they can be calculated with high precision on a computer.

In the case of cover times, the associated Gaussian process of Ding, Lee, Peres also has many symmetries. Formally, it is a stationary Gaussian process. In this case, more simple arguments can be used to calculate the maximum of the Gaussian process. Ultimately, this allows us to compute the cover time as an integral of an important quantity, which we call the volume growth function. We can compute this function with our formulas for hitting times. We believe the volume growth function is the right quantity to understand in order to understand cover times, and is computed on a case by case basis. This method is extremely general, but it faces the difficulty of being very hard to compute since the behavior of all hitting times are needed to understand the volume growth function.

## 2 Background

We review the basics of random walks on groups and Fourier transforms on finite groups that we will later use for the following sections.

### 2.1 Random Walks on Groups

These are very basic facts about random walks on groups that are needed for this paper. See [5] for a more in depth discussion.

**Definition 2.1.** Let $G$ be a group. Let $p$ be a probability measure on $G$. A *random walk on a group* $G$ generated by $p$ is a Markov chain with state space $G$ with the following transition probabilities. For $g, h \in G$, the probability of going from $g$ to $h$ is

$$P(g, h) = p(hg^{-1})$$

We note the symmetry in this definition

$$P(g, hg) = P(e, h) = p(h) \tag{2.1}$$

The random walk is irreducible/reversible if it is an irreducible/reversible Markov chain. Note that the random walk is reversible iff the support of $p$ generates $G$.

**Lemma 2.2.** *The stationary distribution $\pi$ of an irreducible random walk on a group $G$ is the uniform distribution.*

*Proof.* For $g \in G$,

$$\sum_{h \in G} P(h, g) \frac{1}{|G|} = \sum_{h \in G} P(h^{-1}g, g) \frac{1}{|G|} = \frac{1}{|G|} \sum_{h \in G} p(h) = \frac{1}{|G|}$$

□

**Definition 2.3.** For a discrete Markov chain $(X_t)_{t \geq 0}$, let

$$\tau_x = \min\{t \geq 0 : X_t = x\}$$

be the *hitting time*. We will also call

$$\mathbf{E}_y[\tau_x]$$

hitting times. For $\tau_x^+ = \min\{t > 0 : X_t = x\}$, we call

$$\mathbf{E}_x[\tau_x^+]$$

the *first return time*.

**Corollary 2.4.** *For an irreducible random walk on a group $G$. The first return time $\mathbf{E}_e[\tau_e^+]$ of a random walk on a group is $|G|$.*

*Proof.* If $\pi$ is the stationary distribution, then

$$\mathbf{E}_e[\tau_e^+] = \frac{1}{\pi(e)}$$

(see [5] section 1.5). □

**Lemma 2.5.** *A random walk on a group $G$ generated by distribution $p$ is reversible iff $p(g) = p(g^{-1})$ for all $g \in G$.*

*Proof.* Let $\pi$ be the stationary distribution. A Markov chain is reversible iff for all $g, h$

$$\pi(g)p(g,h) = \pi(h)p(h,g)$$

By lemma 2.2, this is

$$\frac{1}{|G|}p(hg^{-1}) = \frac{1}{|G|}p(gh^{-1})$$

which is true iff $p(g) = p(g^{-1})$ for all $g \in G$. □

*Remark* 2.6. A reversible random walk on a group $G$ is a random walk on the Cayley graph with edge weights given by $p$. (This is true for random walks that are not reversible for a directed Cayley graph.)

## 2.2 Fourier Transform on Finite Groups

We review the basics of Fourier transforms on finite groups which will be used in the next section. Proofs will be omitted, but can be found in [2].

**Theorem 2.7** (Fourier Transform). *Let $G$ be a finite group and $f : G \to \mathbb{R}_{\geq 0}$. Then we define the Fourier transform $\hat{f}$ as a function that takes representations $\rho$ of $G$ and gives complex matrices according to the formula*

$$\hat{f}(\rho) = \sum_{g \in G} f(g)\rho(g)$$

**Definition 2.8.** Let $f_1, f_2 : G \to \mathbb{R}_{\geq 0}$ and let $G$ be a finite group. Then the *convolution* is

$$(f_1 * f_2)(h) := \sum_{g \in G} f_1(hg^{-1})f_2(g)$$

*Remark* 2.9. We note that for a random walk on a group with transition probability $P$ and initial distribution $Q$. The distribution at time one is $P * Q$.

*Remark* 2.10. Convolution is in non-commutative in general.

As in standard Fourier analysis, the Fourier transform takes convolutions to products.

**Lemma 2.11.**
$$\widehat{f_1 * f_2}(\rho) = \hat{f}_1(\rho)\hat{f}_2(\rho)$$

**Theorem 2.12** (Fourier Inversion). *Let $f : G \to \mathbb{R}_{\geq 0}$ and $G$ a finite group. Let $\rho_1, \ldots, \rho_r$ be the irreducible representations of $G$. Let $d_i$ be the degree of $\rho_i$. Then for $g \in G$,*

$$f(g) = \frac{1}{|G|}\sum_{i=1}^{r} d_i \operatorname{Tr}(\rho_i(g^{-1})\hat{f}(\rho_i))$$

*Remark* 2.13. The Fourier inversion formula shows that $f$ is determined by how $\hat{f}$ acts on the irreducible representations.

**Theorem 2.14.** *Consider a random walk on a group $G$ generated by a probability distribution $p$. Let $P$ be the transition matrix. Then*

$$P = \phi^{-1}D^*\phi$$

*where $D$ is the block diagonal matrix with blocks of the form $\hat{p}(\rho_i)$ and $\phi$ are the change of basis matrices into the Fourier basis.*

See section 3E of [2] for a more detailed statement of the theorem and a proof.

The following is a lemma of representation theory we will use in the next section. We provide it here as to not interrupt the flow of the next section

**Lemma 2.15.** *Let $\rho$ be a non-trivial irreducible representation of a finite group $G$ on a vector space $V$. Then,*

$$\sum_{g \in G} \rho(g) = 0$$

*Proof.* Let $A = \sum_{g \in G} \rho(g)$ be an endomorphism on $V$. Then $\rho(g)A = A = A\rho(g)$, so by Schur's lemma, $A = \lambda I$ is a multiple of the identity. $\operatorname{Tr} A = \sum_{g \in G} \chi(g) = 0$, by the orthogonality relations for characters, so $A = 0$. □

## 3 Hitting Times

We will prove the following formula for the hitting time (definition 2.3) for random walks on a group using Fourier transforms. Our main idea was that the Fourier basis "almost diagonalizes" matrices of a group, and is thus the natural basis to compute the matrix of hitting times.

**Theorem 3.1.** *Let $G$ be an abelian group, and $\rho_1, \ldots, \rho_r$ be the irreducible representations with $\rho_1$ the trivial representation. Let $p$ a probability distribution on $G$, and consider the random walk on $G$ generated by $p$. Then,*

$$\mathbf{E}_e[\tau_g] = \sum_{i=2}^{r} \frac{1 - \rho_i(g^{-1})}{1 - \hat{p}(\rho_i)}$$

We will prove a more general formula for not-necessarily abelian groups. Let $G$ be a group with $|G| = n$. Consider a random walk on $G$ generated by a probability $p$.

**Lemma 3.2.** *The hitting times satisfy the following system of linear equations*

$$\mathbf{E}_e[\tau_e] = 0 \tag{3.1}$$
$$\mathbf{E}_e[\tau_g] = 1 + \sum_{s \in S} \mathbf{E}_s[\tau_g] P(e, s), \quad g \neq e$$

*Proof.* By definition $\mathbf{E}_e[\tau_e] = 0$. Let $(X_t)_{t \geq 0}$ be a random walk on $G$ generated by $p$ s.t. $X_0 = e$ under a probability distribution $\mathbf{P}$. Then

$$\mathbf{E}_e[\tau_g] = \mathbf{E}_\mathbf{P}[\tau_g] = \sum_{s \in G} P(e, s) \mathbf{E}_\mathbf{P}[\tau_g \mid X_1 = s] = \sum_{s \in G} P(e, s) \mathbf{E}_s[\tau_g]$$

□

To analyze equation (3.1), we use the following simplifying notation. For $g, h \in G$ and $e$ the identity element, the hitting times are

$$h(g) := \mathbf{E}_e[\tau_g] = \mathbf{E}_h[\tau_{gh}]$$

Thus, we view these two quantities as vectors instead of matrices.

Writing equation (3.1) in the above notation we get

$$h(e) = 0 \tag{3.2}$$
$$h(g) = 1 + \sum_{s \in S} h(gs^{-1}) p(s), \quad g \neq e$$

This system of equations looks very similar to a convolution. In order to make this true, we rewrite the equation $h(e) = 0$. We calculate

$$1 + \sum_{s \in S} h(es^{-1})p(s)$$

which is $\mathbf{E}_s[\tau_e^+]$ the first return time. By lemma 2.4,

$$\mathbf{E}_e[\tau_e^+] = \frac{1}{\pi(e)} = n$$

Now we can write equations (3.2) as

$$h(e) = 1 - n + \sum_{s \in S} h(es^{-1})p(s) \tag{3.3}$$

$$h(g) = 1 + \sum_{s \in S} h(gs^{-1})p(s), \quad g \neq e$$

We take the constant terms in the above system, and create a function

$$k(g) = \begin{cases} 1 - n, & g = e \\ 1, & \text{otherwise} \end{cases}$$

The equations (3.3) become

$$h(g) = k(g) + (h * p)(g) \tag{3.4}$$

Let $\rho_1, \ldots, \rho_r$ be the irreducible representations of $G$. Now we may take the Fourier transform to get

$$\hat{h}(\rho_i) = \hat{k}(\rho_i) + \hat{h}(\rho_i)\hat{p}(\rho_i) \tag{3.5}$$

or equivalently

$$\hat{h}(\rho_i)(I - \hat{p}(\rho_i)) = \hat{k}(\rho_i) \tag{3.6}$$

where $I$ is the identity matrix of appropriate dimension.

We compute

$$1 - \hat{p}(\rho_1) = 1 - \sum_{g \in G} p(g)\rho_1(g) = 1 - \sum_{g \in G} p(g) = 0$$

and thus we cannot solve for $\hat{h}(\rho_1)$ with this equation. However, we can solve for $\hat{h}(\rho_i)$ with this equation for all $i \neq 1$.

**Lemma 3.3.** *Consider the random walk on a group $G$ generated by the distribution $p$. If the random walk is irreducible, then $I - \hat{p}(\rho_i)$ is invertible for all nontrivial irreducible representations.*

*Proof.* Let $P$ be the matrix of transition probabilities. The solutions $v$ to $v(I - P) = 0$ are multiples of the stationary distribution. This distribution is unique, since the Markov chain is irreducible. Thus, the nullspace of $I - P$ is 1, so the matrix has rank $n - 1$.

By lemma 2.14, $I - P = \phi^{-1}D^*\phi$, where $D$ is the block diagonal matrix with $I - \hat{p}(\rho_i)$ along the diagonal, and where $\phi$ is a nonsingular matrix. Since $I - P$ has rank $n - 1$, $D$ must have rank $n - 1$. Since $1 - \hat{p}(\rho_1) = 0$, the rest of the blocks must be nonsingular. $\square$

Thus from equation (3.6), we have

$$\hat{h}(\rho_i) = \hat{k}(\rho_i)(I - \hat{p}(\rho_i))^{-1}, \quad i \neq 1 \tag{3.7}$$

**Lemma 3.4.** *Let $\rho_i$ be a non-trivial representation. Then,*
$$\hat{k}(\rho_i) = -nI$$

*Proof.* By lemma 2.15,
$$\hat{k}(\rho_i) = \sum_{g \in G} k(g)\rho_i(g) = \sum \rho_i(g) - n\rho_i(e) = -nI$$

□

From equation (3.7), we get the equation
$$\hat{h}(\rho_i) = -n(I - \hat{p}(\rho_i))^{-1}, \quad i \neq 1 \tag{3.8}$$

We may now compute $\hat{h}(\rho_1)$. Using Fourier inversion, we have
$$0 = h(e) = \sum_{i=1}^{r} d_i \operatorname{Tr}(\hat{h}(\rho_i))$$

so
$$\hat{h}(\rho_1) = -\sum_{i=2}^{r} d_i \operatorname{Tr}(\hat{h}(\rho_i))$$

We summarize our results in the following theorem.

**Theorem 3.5.** *Conisder an irreducible random walk on a finite group $G$. Let $h$ be the hitting times*
$$h(g) := \mathbf{E}_e[\tau_g] = \mathbf{E}_h[\tau_{gh}]$$

*We have a formula for the Fourier transform $\hat{h}$, and thus we can compute $h$ with Fourier inversion. Let $\rho_1, \ldots, \rho_r$ be the irreducible representations of $G$ with $\rho_1$ the trivial representation. The formula is*
$$\hat{h}(\rho_i) = -n(I - \hat{p}(\rho_i))^{-1}, \quad i \neq 1$$
$$\hat{h}(\rho_1) = -\sum_{i=2}^{r} d_i \operatorname{Tr}(\hat{h}(\rho_i))$$

**Theorem 3.6.** *Let $G$ be an abelian group. Using the same notation from the above theorem,*
$$h(g) = \sum_{i=2}^{r} \frac{1 - \rho_i(g^{-1})}{1 - \hat{p}(\rho_i)}$$

*Proof.* Abelian groups have one-dimensional representations, so
$$\hat{h}(\rho_1) = -\sum_{i=2}^{r} \hat{h}(\rho_i)$$

where
$$\hat{h}(\rho_i) = -\frac{n}{1 - \hat{p}(\rho_i)}$$

We get the desired equation by plugging into the Fourier inversion formula. □

**Example 3.7** (Cycle). We compute the hitting times for the random walk on $\mathbb{Z}_n$, $n \geq 2$ generated by the probability distribution $p$ assigning probability $1/2$ to $1, n-1 \in \mathbb{Z}_n$. (We also may define a random walk on $\mathbb{Z}_2$ by assigning probability 1 to $1 \in \mathbb{Z}_2$. The formulas we derive will still hold in this case.) This is an irreducible Markov chain. We will show the hitting times

$$h(k) = k(n-k)$$

For another proof of this fact, see [5]. The representations of $\mathbb{Z}_n$ are for $0 \leq j \leq n-1$,

$$\rho_j(k) = e^{2\pi i jk/n}$$

So

$$\hat{p}(\rho_j) = \sum_{k=1}^{n} \rho_j(k)p(k) = \frac{e^{\frac{2\pi i j}{n}} + e^{\frac{-2\pi i j}{n}}}{2} = \cos\frac{2\pi j}{n}$$

By theorem 3.1,

$$h(k) = \sum_{i=1}^{n-1} \frac{1 - e^{-2\pi i jk/n}}{1 - \cos\frac{2\pi j}{n}} = \sum_{i=1}^{n-1} \frac{1 - \cos\frac{-2\pi jk}{n} + i\sin\frac{-2\pi jk}{n}}{1 - \cos\frac{2\pi j}{n}}$$

Since $h(k)$ is real, we may ignore the imaginary components

$$h(k) = \sum_{i=1}^{n-1} \frac{1 - \cos\frac{-2\pi jk}{n}}{1 - \cos\frac{2\pi j}{n}}$$

Now we simplify this formula.

**Lemma 3.8.**

$$\sum_{i=1}^{n-1} \frac{1 - \cos\frac{-2\pi jk}{n}}{1 - \cos\frac{2\pi j}{n}} = k(n-k)$$

*Proof.* We first note for $k$ not a multiple of $n$,

$$\sum_{i=1}^{n-1} e^{2\pi i k/n} = -1$$

so taking the real part of this equation

$$\sum_{i=1}^{n-1} \cos(2\pi k/n) = -1$$

Let $z = e^{\pi i/n}$, then

$$\sum_{i=1}^{n-1} \frac{1 - \cos\frac{-2\pi jk}{n}}{1 - \cos\frac{2\pi j}{n}} = \sum_{i=1}^{n-1} \frac{1 - \frac{z^{2jk} + z^{-2jk}}{2}}{1 - \frac{z^{2j} + z^{-2j}}{2}} = \sum_{i=1}^{n-1} \frac{(z^{jk} - z^{-jk})^2}{(z^j - z^{-j})^2}$$

$$= \sum_{i=1}^{n-1} (z^{j(k-1)} + z^{j(k-3)} + \cdots + z^{-j(k-1)})^2$$

$$= \sum_{i=1}^{n-1} (z^{2j(k-1)} + z^{-2j(k-1)} + 2(z^{2j(k-2)} + z^{-2j(k-2)}) + \cdots + (k-1)(z^{2j} + z^{-2j}) + k)$$

$$= \sum_{i=1}^{n-1} (2\cos[2\pi j(k-1)] + 4\cos[2\pi j(k-2)] + \cdots + 2(k-1)\cos[2\pi j] + k)$$

$$= -2 - 4 - \cdots - 2(k-1) + k(n-1) = k(n-k)$$

□

**Example 3.9** (Hypercube). Consider the random walk on the group $\mathbb{Z}_2^m$ generated by the probability distribution assigning probability $1/m$ to each of the elements $(0,\ldots,0,1,0,\ldots,0) \in \mathbb{Z}_2^m$. We will show the non-trivial hitting times are all of order $2^m$.

The representations for this group are the Walsh functions, which are defined as follows. For an element $x \in \mathbb{Z}_2^m$, we let $x_i$ be the $i$-th component. For each $k \in \mathbb{Z}_2^m$ we define a representation $\rho_k$ on $x \in \mathbb{Z}_2^m$ by
$$\rho_k(x) = (-1)^{\sum_{i=1}^m x_i k_i}$$
We note that for $k = (0,\ldots,0) =: \mathbf{0}$, this is the trivial representation. For $x \in \mathbb{Z}_2^m$, let $|x|$ be the number of nonzero entries. We compute the Fourier transform
$$\hat{p}(\rho_k) = \sum_{x \in \mathbb{Z}_2^m} p(x)\rho_k(x) = \frac{1}{m}\sum_{i=1}^m (-1)^{k_i} = \frac{m - 2|k|}{m}$$

By theorem 3.1, the hitting time is
$$h(x) = \sum_{k \neq \mathbf{0}} \frac{1 - (-1)^{-\sum_{i=1}^m x_i k_i}}{1 - \frac{m-2|k|}{m}} = m \sum_{k \neq \mathbf{0}} \frac{1 - (-1)^{-\sum_{i=1}^m x_i k_i}}{2|k|}$$

Because the hypercube is symmetric up to relabeling of the coordinates, $h(x)$ depends only on $|x|$. Let $|x| = j$. It suffices to calculate the hitting times for $x = (1,\ldots,1,0\ldots,0)$ the element with $j$ copies of 1 followed by $n - j$ copies of 0. When $k$ has an odd number of 1's among its first $j$ coordinates, $(-1)^{-\sum_{i=1}^m x_i k_i} = -1$. Thus let $S$ be the subset of $\mathbb{Z}_2^m$ s.t. that there are an odd number of 1's among the first $j$ components. Then,
$$h(x) = m \sum_{k \neq \mathbf{0}} \frac{1 - (-1)^{-\sum_{i=1}^m x_i k_i}}{2|k|} = m \sum_{k \in S} \frac{1}{|k|}$$

We count the number of $k \in S$. There must be an odd number $i$ of 1's among the first $j$ components of $k$. For a fixed $i$, there are $\binom{j}{i}$ to pick these ones. For a fixed $l$, there are $\binom{m-j}{l}$ ways to choose these ones among the last $m - j$ elements. Such an element of $S$ has $|k| = i + l$ and there are $\binom{j}{i}\binom{m-j}{l}$ of these. Thus, hitting time is
$$h(x) = m \sum_{i \text{ odd}} \sum_{l=0}^{m-j} \binom{j}{i}\binom{m-j}{l}\frac{1}{l+i} \tag{3.9}$$

There are cases in which we can compute the exact hitting time. We will use the following lemma

**Lemma 3.10.**
$$\frac{2^{m+1} - 1}{m+1} = \sum_{i=0}^m \frac{1}{i+1}\binom{m}{i}$$

*Proof.* We have the binomial formula
$$(x+1)^m = \sum_{i=0}^m \binom{m}{i} x^i$$

Taking the integral $\int_0^x$ to both sides, we get
$$\frac{1}{m+1}(x+1)^{m+1} - \frac{1}{m+1} = \sum_{i=0}^m \binom{m}{i}\frac{x^{i+1}}{i+1}$$

We plug in $x = 1$ to get the identity. □

**Proposition 3.11.** *Let $x \in \mathbb{Z}_2^m$.*
*If $|x| = 1$,*
$$h(x) = 2^m - 1$$

*If $|x| = 2$,*
$$h(x) = \frac{m}{m-1}(2^m - 2)$$

*(The $|x| = 1$ case is proven with a different method in [7].)*

*Proof.* For $|x| = 1$, the formula is

$$h(x) = m \sum_{\substack{i \text{ odd}}} \sum_{l=0}^{m-1} \binom{1}{i}\binom{m-1}{l} \frac{1}{l+i}$$

$$= m \sum_{l=0}^{m-1} \binom{m-1}{l} \frac{1}{l+1}$$

$$= 2^m - 1$$

In the last step, We can do a similar calculation for $|x| = 2$.

$$h(x) = m \sum_{\substack{i \text{ odd}}} \sum_{l=0}^{m-2} \binom{2}{i}\binom{m-2}{l} \frac{1}{l+i}$$

$$= 2m \sum_{l=0}^{m-2} \binom{m-2}{l} \frac{1}{l+1}$$

$$= \frac{m}{m-1}(2^m - 2)$$

□

Now we prove sharp bounds for the hitting times. We start with a lemma

**Lemma 3.12.**
$$\sum_{i+l=s} \binom{j}{i}\binom{m-j}{l} = \binom{m}{s}$$

*where the sum is over all pairs of non-negative integers $(i, l)$ s.t. $i + l = s$.*

*Proof.* The right hand side counts the number of ways to choose $s$ marbles out of $m$. We can count this in another way. Color $j$ marbles red and $m - j$ marbles blue. Then the number of ways to choose $s$ marbles is the sum over all $i + l = s$ of choosing $i$ red marbles and $l$ blue marbles. □

**Proposition 3.13.** *There are non-zero constants $c, C$ such that for $x \in \mathbb{Z}_2^m$,*

$$c2^m \leq h(x) \leq C2^m$$

*for all $x \neq (0, 0, \ldots, 0)$.*

*Proof.* First we prove the upper bound.

$$h(x) = m \sum_{i \text{ odd}} \sum_{l=0}^{m-j} \binom{j}{i}\binom{m-j}{l} \frac{1}{l+i}$$

$$\leq m \sum_{i=0}^{j} \sum_{l=0}^{m-j} \binom{j}{i}\binom{m-j}{l} \frac{1}{l+i}$$

$$\leq m \sum_{s=1}^{m} \sum_{i+l=s} \binom{j}{i}\binom{m-j}{l} \frac{1}{s} \qquad \text{rearranging the sum}$$

$$= m \sum_{s=1}^{m} \binom{m}{s} \frac{1}{s} \qquad \text{by lemma 3.12}$$

$$\leq 2m \sum_{s=1}^{m} \binom{m}{s} \frac{1}{s+1}$$

$$= \frac{2m}{m+1}(2^{m+1} - 1) \qquad \text{by lemma 3.10}$$

Now we prove the lower bound.

**Lemma 3.14.** *First we note that for $i < j/2$, $j$ odd,*

$$3 \sum_{i \text{ odd}} \binom{j}{i} \frac{1}{l+i} \geq \sum_{i=0}^{j} \binom{j}{i} \frac{1}{l+i}$$

*for $l \geq 1$ and*

$$3 \sum_{i \text{ odd}} \binom{j}{i} \frac{1}{i} \geq \sum_{i=1}^{j} \binom{j}{i} \frac{1}{i}$$

*Proof.* First we note that for $i < j/2$, $i + l \geq 2$,

$$\binom{j}{i} \frac{1}{l+i} \geq \binom{j}{i-1} \frac{1}{l+i-1} \Leftrightarrow \frac{j-i+1}{i} \geq \frac{i+l}{i+l-1} \Leftrightarrow il + l + i - 2 \geq 0$$

Thus for $l \geq 1$,

$$\sum_{i \text{ odd}, i<j/2} \binom{j}{i} \frac{1}{i+l} \geq \sum_{i \text{ even}, i<j/2} \binom{j}{i} \frac{1}{i+l} \qquad (3.10)$$

Similarly for $i < j/2$,

$$\binom{j}{i} \frac{1}{l+i} \geq \binom{j}{j-i} \frac{1}{l+j-i}$$

so

$$\sum_{i \text{ odd}, i<j/2} \binom{j}{i} \frac{1}{i+l} \geq \sum_{i \text{ even}, i \geq j/2} \binom{j}{i} \frac{1}{i+l} \qquad (3.11)$$

Using equations (3.10) and (3.11), we have for $j$ odd, $l \geq 1$,

$$3 \sum_{i \text{ odd}, i<j/2} \binom{j}{i} \frac{1}{i+l} + \sum_{i \text{ odd}, i>j/2} \binom{j}{i} \frac{1}{i+l} \geq \sum_{i=0}^{j} \binom{j}{i} \frac{1}{i+l}$$

The same argument holds for $l = 0$. □

Thus by the lemma, we compute the lower bound of the hitting time.

$$\begin{aligned}
h(x) &= m \sum_{\substack{i \text{ odd}}} \sum_{l=0}^{m-j} \binom{j}{i}\binom{m-j}{l} \frac{1}{l+i} \\
&\geq \frac{m}{3} \sum_{i=1}^{j} \sum_{l \leq m-j, i+l \geq 1} \binom{j}{i}\binom{m-j}{l} \frac{1}{l+i} & \text{by lemma 3.14} \\
&\geq \frac{m}{3} \sum_{s=1}^{m} \sum_{i+l=s} \binom{j}{i}\binom{m-j}{l} \frac{1}{s} & \text{by lemma 3.12} \\
&\geq \frac{m}{3} \sum_{s=1}^{m} \sum_{i+l=s} \binom{j}{i}\binom{m-j}{l} \frac{1}{s+1} \\
&= \frac{m}{3} \sum_{s=1}^{m} \binom{m}{s} \frac{1}{s+1} \\
&= \frac{m}{3} \frac{2^{m+1}-1}{m+1} - \frac{m}{3} & \text{by lemma 3.10}
\end{aligned}$$

$\square$

## 4  Cover Times

In this section we would like to use our computations of hitting times to compute cover times defined as follows

$$t_{cov} := \max_{h \in G} \mathbf{E}_h[\max_{g \in G} \tau_g] = \mathbf{E}_e \max_{g \in G}[\tau_g]$$

We will describe our asymptotic bounds with the following notation.

**Definition 4.1.** We say

$$f(n) \lesssim g(n)$$

if there is a universal constant $C$ s.t. $f(n) \leq Cg(n)$. as similarly for $\gtrsim$. If $f(n) \lesssim g(n)$ and $f(n) \gtrsim g(n)$, then we write

$$f(n) \asymp g(n)$$

We will use some results of Ding, Lee, Peres [3] to compute cover times, which we summarize here. Consider a random walk on a group $G$ generated by a distribution $p$. We define the follow Gaussian process $(\eta_g)_{g \in G}$ associated to th random walk, based on the commute time $\kappa(g,h) := \mathbf{E}_g[\tau_h] + \mathbf{E}_h[\tau_g]$. It is the centered Gaussian process with $\eta_e = 0$, and

$$\mathbf{E}[|\eta_g - \eta_h|^2] = \kappa(g,h)$$

These statistics fully determine the Gaussian process.

**Theorem 4.2** (Ding, Lee, Peres). *For a random walk on a group $G$ with vertex set $V$,*

$$t_{cov} \asymp \left( \mathbf{E} \max_{v \in V} \eta_v \right)^2$$

*for the associated Gaussian processes $\eta$ defined above.*

*Remark* 4.3. The results of Ding, Lee, Peres are presented in the more general setting of electrical networks, which we will not use in the result of the paper. For an electrical network, let $c_{gh}$ be the conductance from $g$ to $h$. And $c_g = \sum_{h \in G} c_{gh}$. The total conductance is $\mathcal{C} := \sum_{g,h \in G} c_{gh}$. We define the transition probabilities by $P(g,h) = \dfrac{c_{gh}}{c_g}$. Ding, Lee, Peres define $\mathbf{E}[|\eta_g - \eta_h|^2] = R_{eff} = \dfrac{\kappa(g,h)}{\mathcal{C}}$ (see [2]). Then,

$$t_{cov} \asymp \mathcal{C} \left( \mathbf{E} \max_{v \in V} \eta_v \right)^2$$

Scaling all the conductances, does not change the transition probabilities of the random walk, so $\mathcal{C}$ changes but $t_{cov}$ remains the same. We can also see this in the formula in theorem 4.12 below.

In our case, we let the conductances be $c(g,h) = \dfrac{1}{|G|} p(hg^{-1})$. The total conductance is $\mathcal{C} := \sum_{g,h \in G} c_{gh} = 1$ for this network.

The case of a random walk on a group is much simpler than a general random walk on a graph. This will simplify calculations considerably.

**Definition 4.4.** For a Gaussian process $(\eta_t)_{t \in T}$, $d(t,s) = \mathbf{E}[|\eta_t - \eta_s|^2]^{1/2}$ is a distance on $T$ induced by this process.

In our setting of a random walk on a group with associated Gaussian process, we may check this distance is a metric. Commute time satisfies the triangle inequality because for all $s,t,u \in G$,

$$\mathbf{E}_t[\tau_s] \leq \mathbf{E}_t[\tau_u] + \mathbf{E}_u[\tau_s]$$

Thus

$$d(t,s) = \mathbf{E}[|\eta_t - \eta_s|^2]^{1/2} = \sqrt{\kappa(t,s)} \leq \sqrt{\kappa(t,u) + \kappa(u,s)} \leq \sqrt{\kappa(t,u)} + \sqrt{\kappa(u,s)} = d(t,u) + d(u,s)$$

**Definition 4.5.** A Gaussian process $(\eta_t)_{t \in T}$ is *stationary*, if there is a transitive group action $G \curvearrowright T$ s.t. for $g \in G$,

$$d(g \cdot t, g \cdot s) = d(t,s)$$

**Lemma 4.6.** *Let $G$ be a random walk on a group. The associated Gaussian process is stationary. Namely, the associated metric $d(g,h)$ is fixed by the action of right multiplication by an element of $G$.*

*Proof.* Recall the Gaussian process is defined by for $t,s \in G$,

$$d(t,s) = \mathbf{E}[|\eta_t - \eta_s|^2]^{1/2} = \sqrt{\kappa(t,s)}$$

The group $G$ has a transitive action on itself by right multiplication which fixes $\kappa(t,s)$ since for $g \in G$,

$$\kappa(g \cdot t, g \cdot s) = \mathbf{E}_{tg}[\tau_{sg}] + \mathbf{E}_{sg}[\tau_{tg}] = \mathbf{E}_t[\tau_s] + \mathbf{E}_s[\tau_t] = \kappa(t,s)$$

Thus $d(t,s)$ is fixed by this action. □

**Definition 4.7.** Let $d$ be any metric on a set $T$. Define $N(T,d,\epsilon)$ to be the maximum number of points $t_1,\ldots,t_N \in T$ s.t. $d(t_i,t_j) > \epsilon$ for all $i \neq j$.

**Theorem 4.8** (Fernique)**.** *Let $(\eta_t)_{t \in T}$ be a stationary Gaussian process with finite index set. Then*

$$\mathbf{E}\left[\max_{t \in T} X_t\right] \asymp \int_0^\infty \sqrt{\log N(T,d,\epsilon)} d\epsilon$$

*(See [4] for a proof.)*

Consider a random walk on a group $G$. Combining theorems 4.2 and 4.8, we have

$$t_{cov} \asymp \left(\int_0^\infty \sqrt{\log N(G,d,\epsilon)}d\epsilon\right)^2 \tag{4.1}$$

where $d(g,h) = \sqrt{\kappa(s,t)}$. Now we would like to understand the quantity $N(G,d,\epsilon)$.

**Lemma 4.9.** *Let $G$ be a group with metric $d$. Let $B_g(\epsilon) = \{h \in G : d(g,h) \leq \epsilon\}$ be the closed ball of radius $\epsilon$ around $g$. Then $|B_g(\epsilon)|$ is independent of $g$. We define the quantity $V_d(\epsilon) := |B_g(\epsilon)|$ the volume growth function.*

*Proof.* We will show $|B_g(\epsilon)| = |B_h(\epsilon)|$ for any $g, h \in G$. Let $k \in B_g(\epsilon)$. Then by lemma 4.6,

$$d(g,k) = d(h, kg^{-1}h) \leq \epsilon$$

Right multiplication by $g^{-1}h$ is injective, so it maps $B_g(\epsilon)$ into $B_h(\epsilon)$. Similarly, right multiplication by $h^{-1}g$ maps $B_h(\epsilon)$ into $B_g(\epsilon)$. □

*Remark* 4.10. The term volume growth function already exists in the literature [8]. Here we use the term in a more general context than the existing usage. [8] uses the term volume growth function to refer to the particular case where the distance $d$ is words size. We will be studying a different distance on the group.

**Lemma 4.11.** *Let $G$ be a group with metric $d$,*

$$N(G,d,\epsilon)V_d(\epsilon/3) \leq |G|$$

*and*

$$N(G,d,\epsilon)V_d(\epsilon) \geq |G|$$

*Proof.* First we show the upper bound. Let $S \subset G$ be a set with $N(G,d,\epsilon)$ elements s.t. $d(s,t) \geq \epsilon$ for $s,t \in S$. Then the balls $B_s(\epsilon/3)$ for $s \in S$ are disjoint.

Now we show the lower bound. Let $S$ be as above. We claim the $B_s(\epsilon)$ cover $G$. If this is true, the inequality holds. Assume by contradiction $g \in G$ is not in $B_s(\epsilon)$ for any $s \in S$, then $d(g,s) > \epsilon$ for all $s \in S$, so $S' = S \cup \{g\}$ satisfy $d(s,t) > \epsilon$ for all $s,t \in S'$ and $|S'| > N(G,d,\epsilon)$. However, this contradicts the definition of $N(G,d,\epsilon)$. □

Now, we can understand the cover time in terms of the volume growth function.

**Theorem 4.12.**
$$t_{cov} \asymp \left(\int_0^\infty \sqrt{\log \frac{|G|}{V_d(\epsilon)}}d\epsilon\right)^2$$

*Proof.* From equation (4.1) and lemma 4.11,

$$t_{cov} \asymp \left(\int_0^\infty \sqrt{\log N(G,d,\epsilon)}d\epsilon\right)^2 \gtrsim \left(\int_0^\infty \sqrt{\log \frac{|G|}{V_d(\epsilon)}}d\epsilon\right)^2$$

Similarly, we can get the upper bound

$$t_{cov} \asymp \left(\int_0^\infty \sqrt{\log N(G,d,\epsilon)}d\epsilon\right)^2 \lesssim \left(\int_0^\infty \sqrt{\log \frac{|G|}{V_d(\epsilon/3)}}d\epsilon\right)^2 = \left(3\int_0^\infty \sqrt{\log \frac{|G|}{V_d(\epsilon)}}d\epsilon\right)^2$$

□

**Corollary 4.13** (Matthew's Bounds [6]). *Let $M = \max_{g_1,g_2 \in G} \kappa(g_1, g_2) = \max_{g \in G} \kappa(e, g)$, then*

$$M \lesssim t_{cov} \lesssim M \log |G|$$

*We note $M \asymp \max_{g \in G} \mathbf{E}_e[\tau_g]$ the maximum hitting time.*

*Proof.* The maximum distance $d(s,t) = \sqrt{M}$ for $s, t \in T$. Thus, for $\epsilon \geq \sqrt{M}$, $V_d(\epsilon) = |G|$, so

$$\int_0^\infty \log \frac{|G|}{V_d(\epsilon)} d\epsilon = \int_0^{\sqrt{M}} \log \frac{|G|}{V_d(\epsilon)} d\epsilon$$

Since $V_d(\epsilon) \geq 1$,

$$t_{cov} \lesssim \left(\int_0^{\sqrt{M}} \sqrt{\log \frac{|G|}{V_d(\epsilon)}} d\epsilon \right)^2 \leq \left(\int_0^{\sqrt{M}} \sqrt{\log |G|} d\epsilon \right)^2 = M \log |G|$$

For the lower bound, let $s, t$ be s.t. $d(s, t) = M$. Then $B_s(M/3)$ and $B_t(M/3)$ are disjoint, so $V_d(M/3) \leq |G|/2$.

$$t_{cov} \asymp \left(\int_0^\infty \sqrt{\log \frac{|G|}{V_d(\epsilon)}} d\epsilon \right)^2 \geq \left(\int_0^{\sqrt{M}/3} \sqrt{\log \frac{|G|}{V_d(\epsilon)}} d\epsilon \right)^2 \geq \left(\int_0^{\sqrt{M}/3} \sqrt{\log 2} d\epsilon \right)^2 \gtrsim M$$

□

The following examples show that both the upper and lower bounds may be achieved.

**Example 4.14** (Cycle). We consider the random walk $\mathbb{Z}_n$ generated by the distribution that assigns probability $1/2$ to each of $\pm 1$. We will prove

$$t_{cov} \asymp n^2$$

This is the lower bound of the Matthews bounds, so we must only prove the upper bound. Using example 3.7, $d(0, k) = \sqrt{2k(n-k)}$ for $k \in \mathbb{Z}_n$. Thus for $k < n/2$,

$$\sqrt{2k(n-k)} \leq \epsilon < \sqrt{2(k+1)(n-k-1)} \tag{4.2}$$

implies $V_d(\epsilon) = 2k+1$. The maximum value of $d$ is achieved when $k = \lfloor n/2 \rfloor$. Thus for $\epsilon \geq \sqrt{n^2/2} = n/\sqrt{2}$, $V_d(\epsilon) = n$, so $\log \frac{n}{V(\epsilon)} = 0$. Thus,

$$t_{cov} \asymp \left(\int_0^\infty \sqrt{\log \frac{|G|}{V_d(\epsilon)}} d\epsilon \right)^2 = \left(\int_0^{3\sqrt{n}} \sqrt{\log \frac{|G|}{V_d(\epsilon)}} d\epsilon + \int_{3\sqrt{n}}^{n/\sqrt{2}} \sqrt{\log \frac{|G|}{V_d(\epsilon)}} d\epsilon \right)^2$$

We bound the first summand

$$\int_0^{3\sqrt{n}} \sqrt{\log \frac{|G|}{V_d(\epsilon)}} d\epsilon \leq \int_0^{3\sqrt{n}} \sqrt{\log n} d\epsilon = 3\sqrt{n \log n}$$

Thus, it suffices to show

$$\int_{3\sqrt{n}}^{n/\sqrt{2}} \sqrt{\log \frac{|G|}{V_d(\epsilon)}} d\epsilon \lesssim n$$

Rearranging equation (4.2), we have

$$2(k+1) \geq \frac{\epsilon^2}{n-k-1} \geq \frac{\epsilon^2}{3n}$$

Thus $V_d(\epsilon) \geq \frac{\epsilon^2}{3n} - 1 \geq \frac{\epsilon^2}{6n}$ for $3\sqrt{n} \leq \epsilon \leq \sqrt{2\lfloor n/2 \rfloor^2}$, and we can manually check this inequality for $\sqrt{2\lfloor n/2 \rfloor^2} \leq \epsilon \leq n/\sqrt{2}$. Thus,

$$\int_{3\sqrt{n}}^{n/\sqrt{2}} \sqrt{\log \frac{|G|}{V_d(\epsilon)}} d\epsilon \leq \int_{3\sqrt{n}}^{n/\sqrt{2}} \sqrt{\log \frac{6n^2}{\epsilon^2}} d\epsilon \leq \int_{3\sqrt{n}}^{n/\sqrt{2}} \log \frac{6n^2}{\epsilon^2} d\epsilon$$

$$= (n/\sqrt{2} - 3\sqrt{n})\log 6 + 2(n/\sqrt{2} - 3\sqrt{n})\log n - 2(x\log x - x)\Big|_{3\sqrt{n}}^{n/\sqrt{2}}$$

$$= n\frac{\log 6}{\sqrt{2}} + \sqrt{2}n\log n - \sqrt{2}n\log(n/\sqrt{2}) + \sqrt{2}n + o(n) \asymp n$$

Thus,
$$t_{cov} \asymp n^2$$

The cover times of our other examples can be computed with the following proposition.

**Proposition 4.15.** *Let $f$ be a fixed function. Consider a random walk on a group $G$ with $|G| = N$ s.t. $h(g) \asymp f(N)$ for all $g \neq e$. Then,*
$$t_{cov} \asymp f(N) \log N$$

*Proof.* The maximum hitting time is $\asymp f(n)$, so this is the Matthews upper bound. It suffices to prove a lower bound.

The minimum distnace $\alpha := \min_{g \neq e} d(e, g) \asymp \sqrt{f(n)}$. Then for $\epsilon < \alpha$, $V_d(\epsilon) = 1$. Thus,

$$t_{cov} \asymp \left( \int_0^\infty \sqrt{\log \frac{|G|}{V_d(\epsilon)}} d\epsilon \right)^2 \geq \left( \int_0^\alpha \sqrt{\log \frac{|G|}{V_d(\epsilon)}} d\epsilon \right)^2 = \left( \int_0^\alpha \sqrt{\log |G|} d\epsilon \right)^2 = \alpha^2 \log |G| \asymp f(n) \log n$$

$\square$

**Example 4.16** (Hypercube). From proposition 3.13, $\mathbf{E}_e[\tau_g] \asymp 2^m$ for all $g \neq e$, so by the above proposition the cover time of the hypercube $\mathbb{Z}_2^m$ is

$$t_{cov} \asymp m 2^m$$

**Example 4.17** (Finite-Dimensional Torus). We consider a random walk on a torus $\mathbb{Z}_n^m$. For $n = 2$, we have the hypercube discussed in example 3.9. For $n > 2$, this random walk is generated by the distribution $p$ that assigns probability $1/2m$ to each element $\{0, \ldots, 0, \pm 1, 0 \ldots, 0\}$. We consider the case when $m$ is fixed and $n$ grows. The hitting times for this random walk are computed in section 10.4 of Levin, Peres, Wilmer [5] with electrical network methods. They prove for $m \geq 3$,

$$\mathbf{E}_e[\tau_x] \asymp n^m$$

Thus by proposition 4.15,
$$t_{cov} \asymp n^m \log n$$

# 5 Further Directions

Our method gives us explicit formula for hitting times for random walks on any finite groups whose representations are understood, which includes all finite abelian groups. However, we were only able to analyze and understand these formula for a few simple cases. A generalization of the cases considered would be the torus $\mathbb{Z}_n^m$, letting both $n$ and $m$ grow. This specializes to both the examples of the hypercube (example 3.9) and the finite-dimensional torus (example 4.17). In the cases in this paper, it seems that for large $m$, the hitting times are all of the same order, the order of the group. If this were the case, proposition 4.15, gives us that

$$t_{cov} \asymp mn^m \log n$$

for these cases. This conjecture can be generalized to random walks on a group $G$ generated by a probability distribution $p$. In the above case, we had that the support of $p$ was large enough, and that $p$ was uniform on the support. In these cases, will we also have the phenomenon that the hitting times all are on the order of $|G|$? In these cases where the support of $p$ is large, intuitively it is easy for the random walk to get "lost" even if it starts very close to its destination, which is reason to believe the hitting times would all be of the same order.

A different direction may be the following. The representation theory methods we developed were for the expected values of specific hitting times $\tau_g = \{\min t \geq 0 : X_t = g\}$. Thus, techniques may also be generalized to certain settings for hitting times of the form $\tau_A = \{\min t \geq 0 : X_t \in A\}$ for appropriately chosen sets such as cosets. There is also freedom to choose the initial distribution $\mu$ so that

$$\mathbf{E}_\mu[\tau_A]$$

may be easier to compute and understand.

*I pledge my honor that I have not violated the honor code on this assignment.* -Christopher Zhang